\documentclass[11pt]{amsart}
\usepackage{amsmath, amsthm, latexsym, amssymb, amsfonts, epsfig, epsf}
\usepackage[dvips]{color}

\newenvironment{pf}{\proof[\proofname]}{\endproof}
\theoremstyle{plain}
\newtheorem{Th}{Theorem}[section]
\newtheorem{Cor}[Th]{Corollary}
\newtheorem{Prop}[Th]{Proposition}
\newtheorem{Lemma}[Th]{Lemma}
\numberwithin{equation}{section} \theoremstyle{definition}

\newtheorem{Ex}{Example}
\newtheorem{Prob}{Problem} 

\newtheorem{Def}[Th]{Definition}

\newcommand{\R}{\mathbb R}
\newcommand{\N}{\mathbb N}

\newcommand{\cD}{{\mathcal D}}

\newcommand{\tdet}{{\rm tdet\,}}
\newcommand{\tropdet}{{\rm tropdet\,}}

\newcommand{\rl}[1]{Lemma~\ref{L:#1}}
\newcommand{\rp}[1]{Proposition~\ref{P:#1}}
\newcommand{\rpr}[1]{Problem~\ref{Pr:#1}}

\newcommand{\rex}[1]{Example~\ref{ex:#1}}
\newcommand{\re}[1]{(\ref{e:#1})}
\newcommand{\rc}[1]{Corollary~\ref{C:#1}}
\newcommand{\rt}[1] {Theorem~\ref{T:#1}}

\begin{document}

\author{Thomas Dinitz} 
\author{Matthew Hartman} 
\author{Jenya Soprunova}
\address{Department of Mathematical Sciences,
Kent State University, Kent, OH 44242}
\email{soprunova@math.kent.edu}

\title{Tropical Determinant of Integer Doubly-Stochastic Matrices}

\begin{abstract}  
Let ${\cD }(m,n)$ be the set of all the integer points in the $m$-dilate of the Birkhoff polytope of doubly-stochastic $n\times n$ matrices.  In this paper we find the sharp upper bound on the tropical determinant over the set ${\cD}(m,n)$. We define a version of the tropical determinant where the maximum over all the transversals in a matrix is replaced with the minimum and then find the sharp lower bound on thus defined tropical determinant over ${\cD}(m,n)$.
\end{abstract}
\maketitle

\section*{Introduction}
We start with the following problem. Consider usual Rubik's cube with 9 square stickers on each side colored in one of six colors.  We want to solve Rubik's cube by peeling off the stickers and replacing them so that each of the faces has all stickers of one color. 
Doing this, we will ignore the structure of the Rubik's cube. For example, if initially the (solved) cube had the blue and the green squares on opposite faces, after removing and replacing the stickers
we may end up having blue and green faces adjacent to each other.  Here is our first problem.

\begin{Prob}\label{Pr:rubik}
How many stickers we would need to peel off and replace in the worst case scenario?
\end{Prob}
More generally, assume that we have $n$ pails with $m$ balls in each.  Each ball is colored in one of $n$ colors and we have $m$ balls of each color. Same question: 
\begin{Prob}\label{Pr:pails}
How many balls do we need to move from one pail to another in the worst case scenario so that the balls are sorted by color?
\end{Prob}
Consider an $n\times n$ matrix  $A$  where the rows represent the colors of the stickers (or balls) and the columns represent  the  faces of the cube (or pails). In the $(i,j)$th position in $A$ we record the number of  stickers (balls) of color $i$ on  face (in pail) $j$.  All the entries of matrix $A$
are nonnegative integers and the row and column sums of $A$ are  equal to $m$. We call such matrices {\it integer doubly-stochastic}. 

We would  like to assign each face a color so that the overall number of stickers that we need to move is the smallest possible. In other words, we would like to find a transversal of $A$ with 
the largest possible sum of entries. That is, \rpr{pails} reformulates as 

\begin{Prob}\label{Pr:tdet}
Given positive integers  $m$ and $n$, find  the sharp  lower bound $L(m,n)$ on
$${\rm max}_{\sigma\in S_n}\{a_{1\sigma(1)}+a_{2\sigma(2)}+\cdots+a_{n\sigma(n)}\}
$$
over the set of all integer doubly-stochastic $n\times n$ matrices $A=(a_{ij})$ whose row and column sums are equal to $m$. The answer to \rpr{pails}   would then be $mn-L(m,n)$.
\end{Prob}

The quantity above looks very similar to the determinant of a matrix, where multiplication is replaced with addition while the addition is replaced with taking maximum. Even more, it is  almost identical with the definition of the {\it tropical determinant} 
$$\tropdet A={\rm min}_{\sigma\in S_n}\{a_{1\sigma(1)}+a_{2\sigma(2)}+\cdots+a_{n\sigma(n)}\}
$$
which is related to the classical  assignment problem. Consider $n$ workers and $n$ jobs. Let worker $i$ charge $a_{ij}$ dollars for job $j$. We would like to assign the jobs, one for each worker, so that the overall cost is as small as possible. Clearly, the tropical determinant solves  this problem. A polynomial-time algorithm for solving the assignment problem was developed by Harold Kuhn in 1955 \cite{Kuhn}.

It is now  natural to pose a question similar to \rpr{tdet} where we compute the upper bound on the usual tropical determinant.

\begin{Prob}\label{Pr:tropdet}
Given positive integers  $m$ and $n$, find  the sharp  upper bound $U(m,n)$ on
$$\tropdet A= {\rm min}_{\sigma\in S_n}\{a_{1\sigma(1)}+a_{2\sigma(2)}+\cdots+a_{n\sigma(n)}\}
$$
over the set of all integer doubly-stochastic $n\times n$ matrices $A$ whose row and column sums are equal to $m$. 
\end{Prob}

The set of all doubly-stochastic  $n\times n$ matrices forms a convex polytope  in $\R^{n^2}$  \cite{Birk}, an $m$-dilate of the Birkhoff polytope. The tropical determinant defines a piece-wise linear function on  that  polytope. We would like to to minimize that function over the integer points of the polytope, so the question we are interested in is an integer linear-programming problem. 

In this paper we solve \rpr{tdet} and \rpr{tropdet} completely. First three sections of the paper are devoted to \rpr{tdet}. In the last section we solve \rpr{tropdet}, which turns out to be significantly easier than \rpr{tdet}. Our methods are elementary, combinatorial in nature; our main tool is Hall's marriage theorem.
\vspace{-.3cm}
\subsection*{Acknowledgments} We are thankful to Ivan Soprunov for fruitful discussions.
\section{Definitions, Examples, and Easy Cases}

\begin{Def} Let $A=(a_{ij})$ be an $n\times n$ matrix. Define
$$\tropdet A= {\rm min}_{\sigma\in S_n}\{a_{1\sigma(1)}+a_{2\sigma(2)}+\cdots+a_{n\sigma(n)}\}
$$
$$\tdet A  = {\rm max}_{\sigma\in S_n}\{a_{1\sigma(1)}+a_{2\sigma(2)}+\cdots+a_{n\sigma(n)}\}
$$
\end{Def}

We will refer to both of these quantities as the tropical determinant of $A$, which should not cause confusion since throughout the paper we will mostly be dealing with $\tdet A$ except for the last  section of the paper which is devoted to $\tropdet A$.

\begin{Def}  Let $A$ be an $n\times n$ matrix with non-negative integer entries. We say that $A$ is integer doubly stochastic with sums $m$ if the entries in each of the rows and columns  sum up to a fixed $m\in \N$. We will denote the set of all such matrices by $\cD(m,n)$.
\end{Def}

Then \rpr{tdet} reformulates as
\begin{Prob}\label{Pr:prob1} Fix $m,n\in \N$. Find the sharp lower bound $L(m,n)$ on the tropical determinant  $\tdet A$ over the set $\cD(m,n)$.
\end{Prob}

\begin{Ex} Let $n=5$, $m=7$.
\[\tdet\left(
\begin{array}{ccccc}
1&0&2&{\fbox 2}&2\\
0&1&\fbox{2}&2&2\\
2&\fbox{2}&1&1&1\\
\fbox{2}&2&1&1&1\\
2&2&1&1&\fbox{1}\\
\end{array}
\right)=9
\]
We will later show that $L(7,5)=9$, that is, the minimum of the tropical determinant on the set of $5\times 5$ doubly-stochastic integer matrices with sums 7  is attained on this matrix. 
\end{Ex}

One of our tools is Hall's marriage theorem. The theorem in our formulation deals with a block of zeroes in a matrix with the largest sum of dimensions after all possible swaps of columns and rows. We will refer to such a block as the largest block of zeroes.  

\begin{Th}(Philip Hall \cite{Hall})\label{T:Hall} Let $A$ be an $n\times n$ 0-1 matrix. Then there is a transversal in $A$ that consists of all 1's if and only if the largest  block of zeroes in $A$ has sum of dimensions less than or equal to $n$. 
\end{Th}

Here the theorem is formulated in its weakest form and it can be easily proved by induction on $n$.
For our future discussion we will need two of its corollaries.

\begin{Cor}\label{C:marriage} Let $A$ be an $m\times n$ 0-1 matrix. Then there is a transversal that consists of all 1's if and only if the  largest block of zeroes in $A$  has sum of dimensions less than or equal to $\max(m,n)$. 
\end{Cor}

\begin{pf} Let us assume that $m\geq n$. Extend $A$ to a square 0-1 matrix by appending to $A$ $m-n$ columns consisting of all 1's and apply Hall's marriage theorem to the resulting matrix.
\end{pf}

Let $A$ be an $n\times n$ 0-1 matrix and $W$ be  the block  of zeroes in $A$ with the largest sum of dimensions. Then after some row and column swaps $A$ can be written in the form
\[
A=\left(\begin{array}{cc}   
X&Y\\
Z&W
\end{array}\right)
\]

\begin{Cor}\label{C:Marriage}  Each of $Y$ and $Z$ has a transversal that consists of all 1's.
\end{Cor}
\begin{pf} 
Let $W$ be  of size $d_1$ by  $d_2$ with $d_1+d_2=d$ and let the largest block of zeroes in $Y$ be of  size $s_1$ by  $s_2$. We can assume that the block of zeroes is in the lower right corner of $Y,$ right on top of  the zero block $W$. Then the lower right $s_1+d_1$ by $s_2$ block of $A$ consists of all zeroes and hence  $s_1+d_1+s_2\leq d_1+d_2$, so $s_1+s_2\leq d_2$ and by \rc{marriage}  there exists a  transversal in $Y$ that consists of all 1's. Similarly, such a transversal exists in $Z$.
\end{pf}

Here is the first instance where we are going to apply Hall's marriage theorem to integer doubly-stochastic matrices.
\begin{Lemma}\label{L:nonzero} Let $A\in \cD(m,n)$. Then $A$ has a transversal all of whose entries are nonzero.
\end{Lemma}
\begin{pf} Let's rearrange the rows and columns of $A$ so that the block of zeroes with the largest  sum of dimensions  is in the lower right corner of $A$.  That is, $A$ is of the form
\[
\left(\begin{array}{cc}
{ B}&{\ C}\\
{ D}&
\begin{array}{ccc}
0&\cdots&0\\
\vdots & \ddots & \vdots \\
0&\cdots&0
\end{array}
\end{array}\right).
\]
Let us assume that the block of zeroes is of size $r\times s$.  Let $\Sigma_B$ and $\Sigma_C$ denote the sums of all the entries in the matrices $B$ and $C$ correspondingly. Then we have 
\begin{eqnarray*}
\Sigma_C&=&sm\\
\Sigma_B+\Sigma_C&=&(n-r)m
\end{eqnarray*}
Hence $\Sigma_B=(n-r-s)m$, which implies $r+s\leq n$. By Hall's marriage theorem $A$ has a transversal all of whose entries are nonzero. 
\end{pf}

\begin{Prop}\label{P:monotone}
$L(m,n)\leq L(m+1,n)$.
\end{Prop}
\begin{pf} It's enough to check that for each matrix $A\in\cD(m+1,n)$ there is a matrix $A'\in\cD(m,n)$ such that
$$\tdet A'\leq\tdet A.
$$ 
By \rl{nonzero}, $A$ has a transversal all of whose entries are nonzero. Let $A'$ be obtained from $A$ by subtracting 1 from each element on such a transversal. Clearly, $A'\in\cD(m,n)$ and $\tdet A'\leq\tdet A.$
\end{pf}

\begin{Prop}\label{P:q0}
If $1\leq m\leq n$ then $L(m,n)=n$.
\end{Prop}
\begin{pf} Clearly, $L(1,n)=n$. Hence by \rp{monotone} $L(m,n)\geq n$. To show that $L(m,n)$ equals $n$, we need to construct a matrix a in $\cD(m,n)$ with tropical determinant $n$.  Let the first row of this matrix have $m$ 1's in the first $m$ slots and zeroes in the rest. The second row is the shift of the first row by one slot to the right, etc. For example, here is a matrix in $\cD(4,6)$ whose tropical determinant is 6.
 \[\tdet\left(
\begin{array}{cccccc}
1&1&1&1&0&0\\
0&1&1&1&1&0\\
0&0&1&1&1&1\\
1&0&0&1&1&1\\
1&1&0&0&1&1\\
1&1&1&0&0&1
\end{array}
\right)=6
\]
\end{pf}

\begin{Lemma}
$L(m,n)\geq m$
\end{Lemma}
\begin{pf}
If $A\in\cD(m,n)$, the sum of its entries is $mn$. On the other hand, the sum of all the entries of  $A$ is the sum of its entries in $n$ transversals, which is less than or equal to  $n$ times the tropical determinant of $A$. Hence, $\tdet A\geq m$ and $L(m,n)\geq m$.
\end{pf}

\begin{Prop}\label{P:r0}
If $m$ is a multiple of $n$, then $L(m,n)=m$.
\end{Prop}
\begin{pf} By the previous proposition, we only need to find a matrix in $\cD(m,n)$ whose tropical determinant is $m$. Let $m=qn$ for some $q\in N$ and then the matrix is $qI_n$, where $I_n$ is the $n\times n$ matrix that consists of all 1's.
\end{pf}

\section{Not so Easy Cases}

Let now $m=qn+r$, where $0\leq r<n$. We will provide lower bounds  on the tropical determinant of $A\in \cD(m,n)$, first for the case when $r\geq n/2$, and next for the case when $r<n/2$. Our first bound will turn out to be sharp, while the second bound  will be sharp only under the additional assumption that $qr\geq n-2r$.

\begin{Th}\label{T:easy1} Let $m=qn+r,$ where $n/2 \leq r< n$. Then for any matrix $A\in\cD(m,n)$ we have 
$$\tdet A\geq m+(n-r)=n(q+1).
$$
That is, $L(m,n)\geq n(q+1)$.
\end{Th} 
\begin{pf} Let's assume that there exists a matrix $A\in\cD(m,n)$ such that 
$$\tdet A<m+(n-r)=n(q+1).$$  We rearrange rows and columns  of $A$ so that the tropical determinant is equal to the sum of entries on the main diagonal of $A$ and the entries are decreasing along the main diagonal. That is,
we may assume that

 \begin{equation}
 A= \left( \begin{array}{ccccc}
a_1&   &   &   & b_1 \\
   & a_2 &  &   & \vdots\\
   &   & \ddots &  &\vdots \\
  &  &  & a_{n-1} & b_{n-1}\\
c_1 & \dots & \dots & c_{n-1} &a_{n}\\
\end{array} \right)\in\cD(m,n)
\label{e:matrix}
\end{equation}
\vspace{.2cm}
where $a_1\geq\cdots\geq a_n$.  By our assumption,
$$\tdet A= a_1+\cdots+a_n< n(q+1),
$$
and hence 
$a_n\leq q$ as it is the smallest of the $a_i$'s. We next observe that 
$$c_1+b_1\leq a_1+a_n
$$
since otherwise we could switch the first and the last rows of the matrix and get a bigger sum of entries on the main diagonal. Similarly, we get
\begin{eqnarray*}
c_1+b_1&\leq& a_1+a_n\\
c_2+b_2&\leq& a_2+a_n\\
&\cdots&\\
c_{n-1}+b_{n-1}&\leq& a_{n-1}+a_n
\label{e:system}
\end{eqnarray*} 
Summing  up these inequalities we obtain
\begin{equation}\label{e:abg}
m-a_n+m-a_n\leq \tdet A+(n-2)a_n,
\end{equation}
which implies 
$$2m\leq \tdet A +na_n.
$$
Using our assumption $\tdet A<qn+n$ and its consequence $a_n\leq q$ we get 
$$2qn+2r=2m\leq \tdet A+na_n<qn+n+qn,
$$
so  $2r<n$, which contradicts the hypotheses of the theorem.
\end{pf}

\begin{Th}\label{T:easy2} Let $m=qn+r,$ where $0< r\leq n/2$. Then for any matrix $A\in\cD(m,n)$ we have 
$$\tdet A\geq m+r=qn+2r.
$$ 
That is, $L(m,n)\geq m+r$.
\end{Th} 
\begin{pf}
As in the proof of the previous theorem, we assume that there exists a matrix $A\in\cD(m,n)$ such that 
$$\tdet A<m+r=qn+2r.$$
and that  $A$ is of the form \ref{e:matrix} with non-increasing $a_i$'s, and the sum of entries on the main diagonal equal to the tropical determinant.

Then, using the hypothesis that $2r\leq n$ we get 
$$\tdet A=a_1+\cdots+a_n<qn+2r\leq n(q+1)
$$
which implies $a_n\leq q$ as $a_n$ is the smallest among the $a_i$'s.
As before, we have $2m\leq \tdet A+na_n$, which together with $\tdet A<qn+2r$ and  $a_n\leq q$ implies
$$2qn+2r=2m\leq \tdet A+na_n<qn+2r+qn=2qn+2r,
$$
that is, $2qn+2r<2qn+2r$, a contradiction.
\end{pf}

We next show that the bound obtained in \rt{easy1} is sharp.

\begin{Th}  Let $m=qn+r,$ where $n/2 \leq r<n$. Then
$$L(m,n)=m+n-r=n(q+1).
$$
\label{T:sharp1}
\end{Th}
\begin{pf}
\vspace{-.6cm}
We have already shown in \rt{easy1} that $L(m,n)\geq n(q+1)$. Hence we only need to come up with a matrix $A$ in $\cD(m,n)$ whose tropical determinant is $n(q+1)$. 
Our matrix $A$ will consist of four blocks
\[
A=\left(\begin{array}{cc}
A_1&A_2\\
A_3&A_4
\end{array}\right),
\]
where $A_2\in M(r,n-r)$ and $A_3\in M(n-r,r)$ are matrices all of whose entries are equal to $q+1$, and   $A_4\in M(n-r,n-r)$ has all of its entries equal to  $q$ . The upper left corner block $A_1$ is of size $r$ by $r$ and is constructed in the following way: Let the first $2r-n$ entries in the first row be equal to $q+1$. Here we are using the fact that $0\leq 2r-n<r$, which is equivalent to the assumption of the theorem that $n/2\leq r<n$. The second row is then a shift by one slot of the first row, the third is the shift of the second, etc. For example,  for  $r=5$ and $n=7$ we get
 \[A_1=\left(
\begin{array}{ccccc}
q+1&q+1&q+1&q&q\\
q&q+1&q+1&q+1&q\\
q&q&q+1&q+1&q+1\\\
q+1&q&q&q+1&q+1\\
q+1&q+1&q&q&q+1
\end{array}
\right)
\]
Clearly, $A\in\cD(m,n)$ and its tropical determinant is at most $n(q+1)$ as all its entries are less than or equal to $q+1$. By the \rt{easy1},  the  tropical determinant of $A$  is at  least $n(q+1)$ and hence $\tdet A=n(q+1)$.
\end{pf}

We will show that  under some additional assumptions the bound of \rt{easy2} is sharp.
\begin{Th}\label{T:sharp2} Let $m=qn+r,$ where $0\leq r\leq n/2$, and assume that $qr\geq n-2r$. Then $L(m,n)=m+r=qn+2r$.
\end{Th} 
\begin{pf} We need to construct a matrix  $A$ in $\cD(m,n)$ whose tropical determinant is $qn+2r$.  We will start out in the same way as in \rt{sharp1}. The matrix $A$ will be of the form
\[
A=\left(\begin{array}{cc}
A_1&A_2\\
A_3&A_4
\end{array}\right),
\]
where $A_2\in M(r,n-r)$ and $A_3\in M(n-r,r)$ are matrices all of whose entries are equal to $q+1$, and   $A_4\in M(n-r,n-r)$ has all of its entries equal to  $q$ . Now we need the row and column sums of $A_1$ to be equal to $qr+2r-n$. Notice that now $2r-n\leq 0$. We will first make every entry of $A_1$ to be equal to $q$ and will then distribute $2r-n$ between $r$  entries in each row. For this, let us divide $n-2r$ by $r$ with a remainder to get  $n-2r=rl+r'$, where $0\leq r'<r$. Now  subtract $l$ from each entry in $A_1$ and an extra 1 from $r'$ entries in each row and column. We get

\[A_1= \left( \begin{array}{ccccccc}
q-l-1 & \dots & q-l-1 & q-l & \dots & \dots & q-l \\
q-l & q-l-1 & \ddots & q-l-1 & q-l & \ddots & q-l \\
\vdots & q-l & \ddots & \ddots & \ddots & \ddots & \vdots \\
\vdots & \ddots & \ddots & \ddots & \ddots & \ddots & q-l \\
q-l & \ddots & \ddots & \ddots & \ddots & \ddots & q-l-1 \\
q-l-1 & \ddots & \ddots & \ddots & \ddots & \ddots & \vdots \\
\vdots & \ddots & \ddots & \ddots & \ddots & \ddots & \vdots \\
q-l-1 & \dots & q-l-1 & q-l & \dots & \dots & q-l-1\end{array} \right)\]
where the diagonal band of $q-l-1$'s is of width $r'$.  For this construction to work we need the entries of $A$ to be nonnegative. That is, we need,  $q-l-1\geq 0$ for $r'\neq 0$, and $q-l\geq 0$ for $r'=0$.  This is where we use the additional assumption that $qr\geq n-2r$.  If  $r'=0$,  we get  $qr\geq n-2r=rl$, which implies $qr\geq lr$ and $q\geq l$. Next, let $r'\neq 0$. Then $qr\geq n-2r=rl+r'$, which implies $q\geq l+r'/r$, or $q\geq l+1$.
\end{pf}

\begin{Cor}\label{C:sharp3} Let $m=qn+r$ where  $n/3\leq r<n/2$ and $q\neq 0$. Then $L(m,n)=m+r=qn+2r$. 
\end{Cor}
\begin{pf} We have
$$ \frac{n-2r}{r}\leq \frac{3r-2r}{2}=1
$$
and hence the assumption of \rt{sharp2} that $qr\geq n-2r$ is satisfied for all \linebreak $q\geq 1$.
\end{pf}

\begin{Ex}\label{ex:Rubik} Notice that we have solved the Rubik's cube  \rpr{rubik} stated in the introduction.For this problem, we have $n=6$, $m=9$, and $r=3$, so we are under the assumptions of \rt{sharp1}  and   $L(9,6)=m+n-r=12$, so one needs  to replace $mn-L(m,n)=54-12=42$ stickers  in the worst-case scenario represented by the matrix
 \[\left(
\begin{array}{cccccc}
1&1&1&2&2&2\\
1&1&1&2&2&2\\
1&1&1&2&2&2\\
2&2&2&1&1&1\\
2&2&2&1&1&1\\
2&2&2&1&1&1
\end{array}
\right).
\]
\end{Ex}
\section{Hard cases}
Let $m=nq+r$, where $0\leq r<n$. It remains to take care of the situation when $n>2r+rq$ and $r$ and $q$ are  both nonzero. We will adjust the above construction of a matrix $A$ with small tropical determinant to this case.  The matrix 
$A$ will  still consist of four blocks 
\begin{equation}
A=\left(\begin{array}{cc}
A_1&A_2\\
A_3&A_4
\end{array}\right),
\label{e:amatrix}
\end{equation}
where $A_4$ is a square matrix with all its entries equal to $q$, but its size now could be smaller than $(n-r).$ Sub-matrices $A_2$ and $A_3$ will have entries that are equal to $q$ and $q+1$ and the excess (smaller than in the previous construction) will then be distributed over the matrix $A_1$, whose entries will be less than or equal to $q$. 
\begin{Ex}
Let us look at the example where $n=6$, $m=7$, and $q=1.$ If we use the construction of the previous section we would get the matrix

\[\left(
\begin{array}{cccccc}
-3&2&2&2&2&2\\
2&1&1&1&1&1\\
2&1&1&1&1&1\\
2&1&1&1&1&1\\
2&1&1&1&1&1\\
2&1&1&1&1&1
\end{array}
\right)
\]
with a negative entry. Let's make the block of 1's  $3\times 3$ instead of $4\times 4$, then we won't need to put four 2's in the same row. We get

\begin{equation}
A=\left(
\begin{array}{cccccc}
1&0&2&1&2&1\\
0&1&1&2&1&2\\
2&1&1&1&1&1\\
1&2&1&1&1&1\\
2&1&1&1&1&1\\
1&2&1&1&1&1
\end{array}
\right)
\label{e:l76}
\end{equation}
 and the resulting matrix does not have any negative entries. By \rt{sharp2}, $L(7,6)\geq m+r=8$.
 Notice that
  $\tdet A=10$, which is bigger than the bound of~8. We will show  in 
 \rex{ex76} that $L(7,6)=10$. 
  \end{Ex}
 
 \begin{Ex}
 The block of $q$'s in this construction does not have to be square. For example, if $n=5$ and $m=6$ (and hence $q=1$),  we will show that a matrix whose tropical determinant is the smallest possible has a $3\times 4$ block of $1$'s:
 
 \[\left(
\begin{array}{ccccc}
0&1&2&1&2\\
0&2&1&2&1\\
2&1&1&1&1\\
2&1&1&1&1\\
2&1&1&1&1
\end{array}
\right).
\]
 \end{Ex}
 
For general $m,n$ with $n> 2r+rq$ and nonzero $q$ and $r$ (which implies $n> 3$), we  construct $A$ in a similar way. Let $A_1\in M(l_1,l_2)$, $A_2\in M(l_1,n-l_2)$, $A_3\in M(n-l_1,l_2)$, and $A_4\in M(n-l_1,n-l_2)$ with $l_1,l_2\geq r$.  Set the entries of $A_4$ equal to  $q$.   Let the entries of $A_2$ be $q$'s and $q+1$'s with exactly $r$ $q+1$'s in each column and with the  $q+1$'s evenly distributed among the rows so that the numbers of $q+1$'s in different rows differ by no more than 1. This can be easily accomplished by shifting  the $q+1$'s along the 
columns.  We then order the rows in $A_2$ so that the row sums are non-decreasing as we go from top of the matrix to the bottom.
 Matrix $A_3$ is constructed in a similar way. 
 
 Let us  assume for now that we were able to construct $A_1$ and obtained a matrix $A$ in $\cD(m,n)$ with all the entries less than or equal to $q+1$. The  tropical determinant of such a matrix $A$ is less than equal to  $nq+(l_1+l_2)$. Since we are looking for a matrix in $\cD(m,n)$ with the smallest possible tropical determinant, we want to make $l_1+l_2$ as small as possible.

Let $\Sigma_{A_i}$ be the sum of all the entries of the block  $A_i$ for $i=1\dots 4$.\linebreak Then 
$\Sigma_{A_4}=q(n-l_1)(n-l_2)$. Hence $\Sigma_{A_2}=(n-l_2)m-q(n-l_2)(n-l_1)$ and 
$$\Sigma_{A_1}=l_1m-\Sigma_{A_2}=(l_1+l_2)r+l_1l_2q-nr.
$$
 In order for  $\Sigma_{A_1}$ to be non-negative, we need $l_1$ and $l_2$ to satisfy the inequality
$$(l_1+l_2)r+l_1l_2q\geq nr
$$ 
To minimize $\tdet A$, we would like to find $l_1,l_2\geq r$ that satisfy this inequality with $l_1+l_2$ being as small as possible. If the sum of two integers $l_1+l_2$ is fixed, their product $l_1l_2$ is the largest when $l_1=l_2$ or $l_2=l_1+1$, so we can assume that either $l_1=l_2=l$ or $l_1=l+1$, $l_2=l$.  In the first case, we need $l$ to satisfy  $ql^2+2lr-rn\geq 0$ and in the second,  $ql^2+l(2r+q)+r-rn\geq 0$. For each of these cases let's finish the construction of the block $A_1$ of the matrix $A$.

Let us consider the case when $l$ is the smallest positive integer that satisfies 
$$ql^2+2lr-rn\geq 0$$
and $A_1$ is of size $l$ by $l$.  Notice that our assumption  $n> 2r+rq$  implies that  $l$ is greater then or equal to $r$. Let each entry of $A_1$ be $q$, then the overall excess in the first $l$ rows of $A$ would be $r(n-l)-rl=rn-2lr$.  To check that this excess is positive we need to verify that $n\geq 2l$.  If $n$ is even it's enough to see that if the plug in $n/2$ for $l$ into $ql^2+2lr-rn\geq 0$ the inequality is satisfied and by minimality of $l$ we have $l\leq n/2$. If $n$ is odd we plug in $(n-1)/2$ and again (using $n> 2r+rq> 3r$ and $q\geq 1$) conclude that $n\geq 2l$.

We want to distribute this positive excess $rn-2lr$ ievenly (that is, the row sums differ by no more than 1) among the $l$ rows of $A_1$. Let $ql^2+2lr-rn=a\geq 0$. We have
$$\frac{rn-2lr}{l}=\frac{ql^2-a}{l}=ql-\frac{a}{l},
$$
where $a\geq 0$. Hence if $a=0$ the excess in each row is equal to $ql$ and this can be remedied by letting all the entries of $A_1$ be equal to 0. If $a=ql^2$ then there is no excess and $A_1$ consists of all $q$'s.  

In general, we need to construct $A_1$ whose row and column sums are equal to $a$, where $0\leq a\leq ql^2$  and $a$ is distributed as evenly as possible among the rows and columns of $A_1$ so that this distribution agrees with the distribution of the excess in matrices $A_2$ and $A_3$. For this,  first distribute $a$ among the rows making sure that the row sums are non-increasing as we go from the top of the matrix to be bottom. Then distribute each row sum evenly among the entries swapping entries in the same row to make sure that the distribution over the columns is even. Finally, we swap columns of $A_1$ so that the column sums are  are non-increasing as we go from left to right. We have proved

\begin{Prop}\label{P:const1} Let $m=nq+r$, where $0< r<n$. Additionally, assume that $n>2r+rq$ and  $q$ is nonzero.  Let  $l$ be the smallest positive integer that  satisfies $ql^2+2lr-rn\geq 0$. Then there exists a matrix $A\in\cD(m,n)$ such that $\tdet A\leq nq+2l$. That is, $L(m,n)\leq nq+2l$.
\end{Prop} 

We next consider the case when $l_1=l+1$, $l_2=l$ and $ql^2+l(2r+q)+r-rn=a\geq 0$. 
As before, we first let $A_1$ consist of all $q$'s. The excess over first $(l+1)$ rows is then equal to 
 $$(n-l)r-r(l+1)=nr-2lr-r=l^2q+lq-a
 $$
 and if we divide this by the number of the rows of $A_2$ we get
 $$\frac{nr-2lr-r}{l+1}=lq-\frac{a}{l+1}.
 $$
Similarly, the overall excess over first $l$ columns is $(n-l-1)r-rl=nr-2lr-r$, and when we divide this by the number of columns we get
$$\frac{nr-2lr-r}{l}=\frac{l^2q+lq-a}{l}=(l+1)q-\frac{a}{l}.
$$
It can be checked very similarly to the above that the excess $nr-2lr-r$ is nonnegative, that is, that $n\geq 2l+1$ if $n>2$. Hence $0\leq a\leq l(l+1)q$ and we can construct an $l+1$ by $l$ matrix $A_1$ with row and column sums equal to $a$ by distributing $a$ evenly among the entries as before. We have proved

\begin{Prop}\label{P:const2} Let $m=nq+r$, where $0< r<n$. Additionally, assume that $n>2r+rq$ and  $q$ is nonzero. Let  $l$ be the smallest positive integer that  satisfies $ql^2+l(2r+q)-rn\geq 0$. Then there exists a matrix $A\in\cD(m,n)$ such that $\tdet A\leq nq+2l+1$. That is, $L(m,n)\leq nq+2l+1$.
\end{Prop} 

Our last task is to show that the tropical determinant attains its minimum on the matrices that we have constructed.

\begin{Th}\label{T:hard} Let $m=nq+r$, where $0\leq r<n$. Additionally, assume that $n>2r+rq$ and  $q$ and $r$ are nonzero.  Consider two inequalities
\begin{equation}
l^2q+2lr-rn\geq 0
\label{e:eqn1}
\end{equation}
\begin{equation}
l^2q+l(2r+q)+r-rn\geq 0
\label{e:eqn2}
\end{equation}
Let $l$ be the smallest non-negative integer that satisfies at least one of these inequalities. Then
\begin{enumerate}
\item[]Case\,(1) If this $l$ satisfies only \re{eqn2} then $L(m,n)=qn+2l+1$.
\item[]Case\,(2) If this $l$ satisfies both \re{eqn1} and \re{eqn2} then $L(m,n)=qn+2l$.
 \end{enumerate}
\end{Th}

\begin{Ex}\label{ex:ex76} Let us use this theorem to find $L(7,6)$. Here $m=7$, $n=6$, $q=1$, $r=1$. Clearly,  $n>2r+rq$. The two inequalities are $l^2+2l-6\geq 0$ and $l^2+3l-5\geq 0$. The smallest positive $l$ that satisfies at least one of them is $2$ and it happens to satisfy both inequalities.  Hence $L(7,6)=qn+2l=6+4=10$. See \re{l76} above for an example of a matrix in $\cD(7,6)$ that has tropical determinant 10.
\end{Ex}

\begin{pf} 
Notice that the hypotheses  of the theorem imply that $n>2r+2q\geq 2+1=3$, so in what follows we can use the fact that $n>3$.

Let $A\in\cD(m,n)$. We only need to show that $\tdet A\geq qn+2l+1$ in the first case and $\tdet A\geq qn+2l$ in the second case  since we have constructed matrices in Propositions \ref{P:const1} and \ref{P:const2} whose tropical determinants are less than or equal to (and, hence, by this theorem are equal to)  these bounds.

Swap the  rows and columns of $A$ so that the largest (in terms of sum of dimensions) block of elements that are less than or equal to $q$  is in the lower right corner of $A$.  Let  us call this block $W$. Then $A$ is of the form
\[
A=\left(\begin{array}{cc}
X&Y\\
Z&W
\end{array}\right).
\]

Let $X$ be of size $k_1$ by $k_2$.

\begin{Lemma}\label{L:ineq} Under the assumptions of the theorem, in both  cases, we have $r(k_1+k_2)+k_1k_2q\geq rn$. 
\end{Lemma}
\begin{pf} Let $\Sigma_W$ and $\Sigma_Y$  be the sums of all entries in blocks $W$ and  $Y$. Then \linebreak $\Sigma_W\leq q(n-k_1)(n-k_2)$. Hence 
$$\Sigma_Y= (n-k_2)(qn+r)-\Sigma_W\geq (n-k_2)(qn+r)-q(n-k_1)(n-k_2).
$$
On the other hand,  $\Sigma_Y\leq k_1(qn+r)$. Putting these two inequalities together we get
$$(n-k_2)(qn+r)-q(n-k_1)(n-k_2)\leq  k_1(qn+r),
$$
which is equivalent to $r(k_1+k_2)+k_1k_2q\geq rn$.
\end{pf}

Let us assume that the conclusion of the theorem does not hold. That is, we assume in Case (1) that  $\tdet A\leq qn+2l$ and in Case (2) $\tdet A\leq qn+2l-1$. 

\begin{Lemma} Under the assumptions of the proof of the Theorem that we have made so far, $k_1+k_2\leq n$.  
\end{Lemma}

\begin{pf}
If not,  the sum of dimensions  of $W$ is less than $n$ and by Hall's Marriage Theorem \ref{T:Hall} $A$ has a transversal with each of  the entries being at least $q+1$. Hence $\tdet A\geq n(q+1)$. On the other hand,  in Case (2) we are under the assumption that $\tdet A\leq qn+2l-1$, and hence we get $n(q+1)\leq qn+2l-1$  and then $n\leq 2l-1$, which rewrites $l\geq (n+1)/2$, but it is easy to see that $l\leq n/2$. For this, plug in  $n/2$  into \re{eqn1} if $n$ is even and $(n-1)/2$ if $n$ is odd and verify that both of them satisfy the inequality. In Case(1) we get $n(q+1)\leq qn+2l$, that is,  $l\geq n/2$.  Again plugging in $(n-1)/2$ for $l$ into \re{eqn2} if $n$ is odd and $(n-2)/2$ if $n$ is even, we get a contradiction.
\end{pf}

\begin{Lemma}\label{L:kl} Under the assumptions of the proof of the Theorem that we have made so far, we have $k_1+k_2\geq 2l+1$ in Case (1) and  $k_1+k_2\geq 2l$ in Case (2).
\end{Lemma}
\begin{pf}
In Case (2) we have 
 $$(l-1)^2q+(l-1)(2r+q)+r-rn<0,
$$ which rewrites $r(2l-1)+l(l-1)q<rn$. Let us assume that $k_1+k_2\leq 2l-1$. Then $k_1k_2\leq l(l-1)$ and $r(k_1+k_2)+k_1k_2q\leq r(2l-1)+l(l-1)q<rn$ by the above, which contradicts \rl{ineq}.

Similarly, in Case (1)  we have $l^2q+2lr-rn<0$. Let us assume that $k_1+k_2\leq 2l$. Then $k_1k_2\leq l^2$ and $r(k_1+k_2)+k_1k_2q\leq 2lr+l^2q<rn$ by the above, which contradicts \rl{ineq}.
\end{pf}

By \rc{Marriage} there exist  transversals  in each $Y$ and $Z$ whose entries are at least $q+1$. Notice, that since $k_1+k_2\leq n$, the  transversals in $Y$ and $Z$ have respectively $k_1$ and $k_2$ entries.

Consider  transversals in $Y$ and $Z$ that have largest possible sums of entries $t_1$ and $t_2$. Among those, pick one transversal in $Y$ and one in $Z$ such that their union extends to the largest (among transversals of this kind) transversal in $A$. Let $a_1,\dots, a_{k_1}$ and $b_1,\dots, b_{k_2}$ be the entries of those transversals. We have
$$t_1=a_1+\cdots+a_{k_1}\geq k_1(q+1)\  \ \ t_2=b_1+\cdots +b_{k_2}\geq k_2(q+1).
$$

We would like to extend the union of these two transversals to a transversal of $A$.
For this, cross out the rows and columns of $A$ that intersect those largest transversals to get  an $(n-k_1-k_2)$ by $(n-k_1-k_2)$ submatrix $Q$ of $W$. Let us assume first that $Q$ has a transversal of  all $q$'s. 

Then, using \rl{kl},  we get a contradiction in Case (2)
\begin{eqnarray*}
qn+2l-1&\geq& \tdet A \geq  t_1+t_2+(n-k_1-k_2)q\\
&\geq& k_1(q+1)+k_2(q+1)+qn-k_1q-k_2q\\
&=&qn+k_1+k_2\geq qn+2l,
\end{eqnarray*}
as well as  in Case (1) 
\begin{eqnarray*}
qn+2l&\geq& \tdet A \geq  t_1+t_2+(n-k_1-k_2)q\\
&\geq& k_1(q+1)+k_2(q+1)+qn-k_1q-k_2q\\
&=&qn+k_1+k_2\geq qn+2l+1.
\end{eqnarray*}

Next, assume that the largest transversal of $Q$ has an entry that is less than or equal to $q-1$.  We can assume that this largest transversal is the main diagonal of $Q$ and the smallest entry of the largest transversal is in the lower right corner. We have

$$
A= \left( \begin{array}{cccc|cccc|cccc}
&  &  & &        &      &  &c_1           &a_1&&&\\
&  &  & &        &      &  &c_2            &     &a_2&&\\
&  &  & &        &      &  &\vdots         &     &&\ddots&\\
&  &  & &        &      &   &c_{k_1}       &     &&& a_{k_1}\\
\hline
&  &  & &e_1   &                  &  &f_1           & &&&\\
&  &  & &        &e_2 &          &f_2            &     &&&\\
&  &  & &        &      &\ddots  &\vdots         &     &&&\\
d_1&d_2 &\dots  &d_{k_2} &g_1   &g_2 &\dots   &e_{k_3}       &i_1     &i_2&\dots& i_{k_1}\\
\hline
b_1&      &          &            & & &  &h_1           & &&&\\
     &b_2 &          &            & & &  &h_2            &     &&&\\
     &      &\ddots  & &        & & &\vdots         &     &&&\\
     &      &\          &b_{k_2}& &  &  &h_{k_2}       &i_1     &i_2&\dots& i_{k_1}

\end{array} \right)
$$
\vspace{.3cm}

\noindent where $e_{k_3}\leq {q-1}$, $k_3=n-k_1-k_2$, and $\tdet Q=e_1+e_2+\cdots+e_{k_3}$. Let $f=f_1+\cdots+f_{k_3-1}+e_{k_3}$, $g=g_1+\cdots+g_{k_3-1}+e_{k_3}$, $c=c_1+\cdots+c_{k_1}$, $d=d_1+\cdots+d_{k_2}$, $h=h_1+\cdots+h_{k_2}$, and $i=i_1+\cdots+i_{k_1}$. 

By an argument similar to that of \rt{easy1} we get 
\begin{equation}\label{e:W}
f+g\leq k_3e_{k_3}+\tdet Q.
\end{equation}

Next notice that each $c_j\leq a_j$ and $i_j\leq q$. We cannot simultaneously have $c_j=a_j$ and $i_j=q$ since in that case swapping columns we can increase $\tdet Q$ not changing $t=a_1+\cdots+a_{k_1}$. Hence $c_j+i_j\leq a_j+q-1$ and summing these up over $j$ we get $c+i\leq t_1+qk_1-k_1.$ Similarly, $d+h\leq t_2+ qk_2-k_2$.
Since the row and column sums of $A$ are equal to $qn+r$  we get 
$$qn+r=c+f+h,\ \ \ \ qn+r=d+g+i.
$$
Adding up these two equalities  and using inequalities obtained above together with \re{W} we get
\begin{eqnarray*}
2qn+2r&=&(c+i)+(d+h)+(f+g)\\
&\leq& t_1+qk_1-k_1+ t_2+ qk_2-k_2+ k_3e_{k_3}+\tdet Q\\
&\leq& t_1+qk_1-k_1+ t_2+ qk_2-k_2+ k_3(q-1)+\tdet Q\\
&=&t_1+t_2+qn-n+\tdet Q
\end{eqnarray*}
 
From here,
$$\tdet A\geq t_1+t_2+\tdet Q\geq qn+2r+n\geq qn+k_1+k_2 .
$$ 
Next, by \rl{kl}, in Case (1) we have  $k_1+k_2\geq 2l+1$ and hence 
$$\tdet A\geq qn+k_1+k_2\geq qn+2l+1;$$
In Case (2) we have $k_1+k_2\geq 2l$ and hence 
$$\tdet A\geq qn+k_1+k_2\geq qn+2l,$$ 
as claimed.
\end{pf}

We  sum up all the obtained results in
\begin{Th} Let $m=qn+r$ where $0\leq r<n$. Then 
\begin{itemize}
\item If $q=0$  then $L(m,n)=n.$ (\rp{q0})
\item If $r=0$ then $L(m,n)=m.$ (\rp{r0}) 
\item If $n/2\leq r<n$ then $L(m,n)=n(q+1).$ (\rt{sharp1})
\item If $0<r<n/2$ and $n\leq 2r+rq$ (in particular, if $n/3\leq r<n/2$) then $L(m,n)=qn+2r.$  (\rt{sharp2}, \rc{sharp3})
\item If  $n> 2r+rq$ and $r,q\neq 0$ then $L(m,n)=qn+2l+1$ or $qn+2l$. (see the definition of $l$ and details in \rt{hard})
\end{itemize}
\end{Th}

\section{Upper Bound on the  Tropical Determinant}

Let $A=(a_{ij})$ be an $n$ by $n$ matrix. Recall that
$$\tropdet (A)= {\rm min}_{\sigma\in S_n}\{a_{1\sigma(1)}+a_{2\sigma(2)}+\cdots+a_{n\sigma(n)}\}
$$
In this section we solve \rpr{tropdet} which reformulates as
\begin{Prob}\label{Pr:prob2} Fix $m,n\in \N$. Find the sharp upper bound $U(m,n)$ on the tropical determinant $\tropdet A$  over the set $\cD(m,n)$.
\end{Prob} 

This problem turns out to be much easier than our original problem.  The solution is very similar to the ``not so easy'' cases of \rpr{prob1}.

\begin{Th}\label{T:easy3} Let $m=qn+r$ where $0 \leq r\leq n/2$. Then for any matrix $A\in\cD(m,n)$ we have $\tropdet A\geq qn$. That is, $U(m,n)\leq qn$.
\end{Th} 
\begin{pf} Let's assume that there exists a matrix $A\in\cD(m,n)$ such that 
$$\tropdet A>qn.$$  We rearrange rows and columns  of $A$ so that the tropical determinant is equal to the sum of entries on the main diagonal of $A$ and the entries are non-decreasing along the main diagonal. That is, we may assume that $A$ is of the form
\begin{equation}
 A= \left( \begin{array}{ccccc}
a_1 &   &   &   & b_1 \\
   & a_2&  &   & \vdots\\
   &   & \ddots &  &\vdots \\
  &  &  & a_{n-1}& b_{n-1}\\
c_{1} & \dots & \dots & c_{n-1} & a_{n}\\
\end{array} \right)
\label{e:matrix1}
\end{equation}
where $a_1\leq a_2\leq\cdots\leq a_n$ and $\tropdet A=a_1+\cdots+a_n$. Notice that we also have 
By our assumption,
$$\tropdet A= a_1+\cdots+a_n>qn,
$$
which implies $a_n\geq q+1$ as $a_n$ is the largest among the $a_i$'s.
Similarly to the proof of \rt{easy1} we get $2m\geq \tropdet A+na_n$ and hence
$$2qn+2r=2m\geq \tropdet A+na_n>qn+n(q+1),
$$
which implies $2r>n$, and this contradicts the hypotheses of the theorem. 
\end{pf}

\begin{Th}\label{T:easy4} Let $m=qn+r$ where $n/2\leq r< n$. Then for any matrix $A\in\cD(m,n)$ we have $\tropdet A\leq qn+(2r-n)$. That is, $U(m,n)\leq qn+(2r-n)$.
\end{Th} 
\begin{pf}
As in the previous theorem, we assume that $A$ is of the form \ref{e:matrix1} with non-decreasing $a_i$'s and the sum of entries on the main diagonal equal to the tropical determinant.

Let us assume that $\tropdet A>qn+(2r-n)$. Then
$$\tropdet A=a_1+\cdots+a_n>qn+(2r-n)\geq qn,
$$
which implies that $a_n\geq q+1$ as $a_n$ is the largest among the $a_i$'s.
As before we get
$$2m\geq \tropdet A+na_n>qn+(2r-n)+n(q+1)=2qn+2r=2m,
$$
that is, $2m<2m$, a contradiction.
\end{pf}

Finally, we show that these bounds are sharp.

\begin{Th} Let $m=qn+r$ where $0 \leq r<n/2$. Then $U(m,n)=qn$.
\label{T:sharp3}
\end{Th}
\begin{pf}
We have already shown in \rt{easy1} that $U(m,n)\leq nq$. Hence we only need to come up with a matrix $A$ in $\cD(m,n)$ whose tropical determinant is $nq$. 
Our matrix $A$ will consist of four blocks
\[
A=\left(\begin{array}{cc}
A_1&A_2\\
A_3&A_4
\end{array}\right),
\]
where $A_2\in M(n-r,r)$ and $A_3\in M(r,n-r)$ are matrices all of whose entries are equal to $q$, and   $A_4\in M(r,r)$ has all of its entries equal to  $q+1.$ The upper left corner block $A_1$ is of size $n-r$ by $n-r$ and is constructed in the following way: Let the first $r$ entries in the first row be equal to $q+1.$ Here we are using the fact that $r<n-r$, which is equivalent to the assumption of the theorem that $r<n/2$. The second row is then a shift by one slot of the first row, the third is the shift of the second, etc. 

Clearly, the tropical determinant of $A$ is at least $nq$ as all its entries are greater than or equal to $q$. By the \rt{easy3},  the  tropical determinant of $A$  is at  most $nq$ and hence $\tropdet A=nq$.
\end{pf}

\begin{Th}\label{T:sharp4} Let $m=qn+r$ where $n/2\leq r<n$. Then $U(m,n)=qn+2r-n$.
\end{Th} 
\begin{pf} We need to construct a matrix  $A$ in $\cD(m,n)$ whose tropical determinant is $qn+2r$.  As before, let matrix $A$  be of the form
\[
A=\left(\begin{array}{cc}
A_1&A_2\\
A_3&A_4
\end{array}\right),
\]
where $A_2\in M(n-r,r)$ and $A_3\in M(r,n-r)$ are matrices all of whose entries are equal to $q$, and   $A_4\in M(r,r)$ has all of its entries equal to  $q+1.$ It remains to define $A_1$ so that the sums of first $n-r$ rows and columns of $A$ are equal to $qn+r$. Let first all the entries of $A_1$ be equal to $q$ and then distribute $r$ between $n-r$  entries in each row and column of $A_1$. For this, divide $r$ by $n-r$ with a remainder to get  $r=(n-r)q'+r'$, where $0\leq r'<r$.  Notice that $q'\geq 1$ since $r\geq n-r$. Next, let the  entries of $A_4$ be equal to $q+q'$ and then add a band of 1's of width $r'<r$ to $A_1$.  To get  the smallest transversal, pick a transversal of all $q's$ in each $A_2$  and $A_3$ and then extend it to a transversal of $A$. We get
$$\tropdet A=qn+n-2(n-r)=qn+2r-n.
$$
\end{pf}

\end{document}